\documentclass[12pt,a4paper,reqno]{amsart}
\usepackage{a4wide,amsmath,amscd,amssymb,latexsym}

\def\q{\hskip0.17cm}
\def\,{\hskip0.12cm}

\newtheorem{thm}{Theorem}[section]
\newtheorem{cor}[thm]{Corollary}
\newtheorem{lem}[thm]{Lemma}
\newtheorem{prop}[thm]{Proposition}
\theoremstyle{definition}
\newtheorem{Def}[thm]{Definition}
\newtheorem{Exm}[thm]{Example}
\theoremstyle{remark}
\newtheorem{rmk}[thm]{Remark}
\newcommand{\e}{\mathfrak{e}}
\newcommand{\h}{\mathfrak{H}}
\newcommand{\sm}{\left(\smallmatrix}
\newcommand{\esm}{\endsmallmatrix\right)}
\newcommand{\mat}{\begin{pmatrix}}
\newcommand{\emat}{\end{pmatrix}}

\begin{document}
 \title{Traces of singular values and Borcherds products}
 \author{Chang Heon Kim}
 \begin{abstract}
 Let $p$ be a prime for which the congruence group
 $\Gamma_0(p)^*$ is of genus zero, and
 $j_p^*$ be the corresponding Hauptmodul.
  Let $f$ be a nearly holomorphic modular form of
 weight $1/2$ on $\Gamma_0(4p)$ which satisfies some congruence
 condition on its Fourier coefficients. We interpret $f$ as
 a vector valued modular form. Applying Borcherds
 lifting of vector valued modular forms we construct infinite
 products associated to $j_p^*$ and extend Zagier's trace formula
 for singular values of $j_p^*$. Further we investigate the
 twisted traces of sigular values of $j_p^*$ and construct
 Borcherds products related to them.
 \end{abstract}
 \address{
 Mathematisches Institut der Universit\"at Heidelberg,
 Im Neuenheimer Feld 288, 69120 Heidelberg}
 \email{kim@mathi.uni-heidelberg.de}
 \maketitle
 \renewcommand{\thefootnote}%
             {}

\section{Introduction}
 \par
 Let $d$ denote a positive integer congruent to 0 or 3 modulo 4.
 We denote by ${\mathcal Q}_d$ the set of positive definite binary quadratic
 forms $Q=[a,b,c]=aX^2+bXY+cY^2 \, (a,b,c\in \mathbb Z)$ of
 discriminant $-d$, with usual action of the modular group
 $\Gamma=SL_2(\mathbb Z)$. To each $Q\in {\mathcal Q}_d$, we associate its
 unique root $\alpha_Q\in {\mathfrak H}$ (=upper half plane).
 Let $j(\tau)$ ($\tau\in {\mathfrak H}$) be the elliptic modular invariant
 and ${\bf t}(d)$ be the (weighted) trace of a singular modulus of discriminant
 $-d$ ($=\underset{Q\in
 {\mathcal Q}_d/\Gamma}{\sum} \frac{1}{|\bar{\Gamma}_Q|}
  (j(\alpha_Q)-744)$).
  Here $\bar{\Gamma}_Q=\{ \gamma\in \bar{\Gamma}=
  PSL_2(\mathbb{Z}) \mid Q\circ \gamma=Q \}$.
  In addition we set $\textbf{t}(-1)=-1,
  \textbf{t}(0)=2$ and $\textbf{t}(d)=0$ for $d<-1$ or $d\equiv
  1,2$ (mod 4).
 Zagier showed the series
 $\sum_{d\in\mathbb{Z}} {\bf t}(d) q^d$ $
 (q=e^{2 \pi i \tau}, \tau \in \h)$
  is a modular form of weight 3/2 on
 $\Gamma_0(4)=\{ \sm a&b \\ c&d \esm \in SL_2(\mathbb Z), 4 | c
 \}$, holomorphic in $\mathfrak H$ and meromorphic at cusps
 (\cite{Zagier} Theorem 1).
 \par
 Let $\Gamma_0(N)^*$ ($N=1,2,3,\dots$)
 be the group generated by $\Gamma_0(N)$ and all
 Atkin-Lehner involutions $W_e$ for $e||N$. Here
 $W_e$ is represented by a matrix of the form
 $\frac{1}{\sqrt{e}} \sm ex&y\\Nz&ew \esm$ with $x,y,z,w\in
 \mathbb{Z}$ and $xwe-yzN/e=1$.
 There are only finitely
 many values
 of $N$ for which $\Gamma_0(N)^*$ is of genus $0$. In particular,
 if we
 let $\mathfrak{S}$ denote the set of prime values for such $N$,
 then
 $$ \mathfrak{S}=\{2,3,5,7,11,13,17,19,23,29,31,41,47,59,71\}.$$
 Let $j_N^*$ be the corresponding {\em Hauptmodul},
 whose Fourier expansion starts with $q^{-1}+0+a_1 q+a_2 q^2+\cdots$.
 It can be
 described by means of Dedekind eta function or theta functions.
 For example, if $N-1$ divides 24, then the Hauptmodul
 $j_N^*$ is explicitly given by
 $$j_N^*(\tau)
 =\left(\frac{\eta(\tau)}{\eta(N\tau)}\right)^{\frac{24}{N-1}}
 +\frac{24}{N-1}
    +N^\frac{12}{N-1}
    \left(\frac{\eta(N\tau)}{\eta(\tau)}\right)^{\frac{24}{N-1}} .$$
 By a result of Borcherds \cite{B92} the Fourier coefficients of
 $j_N^*$ are related to certain representations of the monster
 simple group.
 \par
 Let $d$ be an integer $\ge 0$ such that $-d$ is
 congruent to a square modulo $4N$.
 We choose an integer $\beta \, (\mbox{mod }
 2N)$ with $\beta^2 \equiv -d \, (\mbox{mod } 4N)$ and consider
 the set $\mathcal{Q}_{d,N,\beta}
 =\{ [a,b,c]\in \mathcal{Q}_d \mid
 a \equiv 0 \, (\mbox{mod } N)$,
 $b \equiv \beta \, (\mbox{mod } 2N)\}$ on which $\Gamma_0(N)$ acts.
 We assume that $d$ is not divisible as a discriminant by the
 square of any prime dividing $N$.
 Then
 the the image of the root $\alpha_Q$ ($Q\in
 \mathcal{Q}_{d,N,\beta}$) in $\Gamma_0(N)\backslash \h$
 corresponds to Heegner point
 and the natural map from
 the quotient $\mathcal{Q}_{d,N,\beta}/\Gamma_0(N)$ to
 $\mathcal{Q}_{d}/\Gamma$ is a bijection. We define a
 trace $\textbf{t}^{(N,\beta)}(d)$ by
 $\textbf{t}^{(N,\beta)}(d)=\sum_{Q\in\mathcal{Q}_{d,N,\beta}/\Gamma_0(N)}
  \frac{1}{|\bar\Gamma_Q|} j_N^* (\alpha_Q)$.
 Since $j_N^*$ is the Hauptmodul for $\Gamma_0(N)^*$, the value of
  $\textbf{t}^{(N,\beta)}(d)$ is independent of
  the choice of $\beta$ and we
 can simply write $\textbf{t}^{(N)}(d)$
 instead of $\textbf{t}^{(N,\beta)}(d)$.
 Zagier  described $\textbf{t}^{(N)}(d)$ as the
 Fourier coefficient of certain Jacobi form of weight 2 and index
 $N$ for $2\le N \le 6$ (see \cite{Zagier} Theorem 8).
 \par
 Let $M_{k+1/2}^{!}(N)$ be the vector space consisting
 of {\em nearly holomorphic modular forms} (holomorphic in $\h$ and
 meromorphic at cusps) of half-integral weight
 $k+1/2$ on $\Gamma_0(4N)$ whose $n$-th Fourier coefficient
 vanishes unless $(-1)^{k}n$ is congruent to a square modulo $4N$.
 According to Borcherds' duality theorem (Theorem 3.1 in
 \cite{B99}), the only obstructions to finding such forms of
 weight 1/2
 are given by holomorphic vector valued modular forms
 of weight $3/2$, which can be
 identified with elements in $J_{2,N}$
 (=the space of holomorphic
 Jacobi forms of weight $2$ and index $N$)
 (see \cite{E-Z} Theorem 5.1).
 Let $N$ be a positive integer for which the dimension of
 $J_{2,N}$ is zero.
 Then
 for every integer $d\ge 0$
 such that $-d$ is congruent to a square
 modulo $4N$, we can find
 a unique modular form $f_{d,N}\in M_{1/2}^{!}(N)$
 having a Fourier expansion of the form
 $$ f_{d,N}=q^{-d}+\sum_{D>0} A(D,d) q^D. $$
 We define a polynomial
 $\mathcal
 H_{d,N}(X)$ by
 $ \underset{Q\in {\mathcal Q}_{d,N,\beta}/\Gamma_0(N)}{\prod}
 (X-j_N^*(\alpha_Q))^{1/|\bar\Gamma_Q|}$.
 In \cite{Kim}, from Zagier's trace formula
 (\cite{Zagier} Theorem 8) we derived the
 following product identity: for $N=2,3,5,6$
 \begin{equation} \label{PRODUCT}
 {\mathcal H}_{d,N}(j_N^*(\tau))=q^{-H(d)}\prod_{u=1}^\infty
 (1-q^u)^{A^*(u^2,d)}
 \end{equation}
 where $H(d)$ ($=\underset{Q\in {\mathcal Q}_d/\Gamma}{\sum}
 \frac{1}{|\bar\Gamma_Q|})$ is
 the \textit{Hurwitz-Kronecker class number} and
 $A^*(D,d)$ is defined to be $2^{s(D,N)} A(D,d)$ with $s(D,N)$
 the number of distinct prime factors dividing $(D,N)$.

 \par
 Given an even lattice $M$ in a real quadratic space
 of signature $(2,n)$, Borcherds
 lifting \cite{B98} gives a multiplicative correspondence between
 vector valued modular forms $F$ of weight $1-n/2$
 with values in
 $\mathbb{C}[M'/M]$ (= the group ring of $M'/M$)
 and meromorphic modular forms on complex
 varieties $(O(2)\times O(n)) \backslash O(2,n) / {\rm Aut}(M,F)$.
 Here $O(2,n)$ is the orthogonal group of $M\bigotimes \mathbb{R}$
 and
 ${\rm Aut}(M,F)$ is the subgroup of orthogonal group of $M$
 ($={\rm Aut}(M)$)
 leaving the form $F$ stable under the
 natural action of ${\rm Aut}(M)$ on $M'/M$.
 In this article we are mainly concerned with $O(2,1)$ case.
 For $O(2,n)$ ($n\ge 2$) case one is referred to
 \cite{Bruinier, Br-Bu}.
 We take $M$ to be the direct sum of 2-dimensional unimodular
 Lorentzian lattice and 1-dimensional lattice generated by a
 vector of norm $2N$, as in
 Example 5.1 of \cite{B99}.
 If $N=1$, each element in $M^!_{1/2}(1)$ has a natural
 interpretation as a vector valued modular form
 (\cite{B95} Lemma 14.2). Let $N$ be a prime $p$
 for which $\Gamma_0(p)^*$ is of genus zero, i.e. $p\in \mathfrak{S}$.
 We will interpret a scalar valued modular form $f\in M^!_{1/2}(p)$
 as a vector valued modular
 form $F$ with values in
 $\mathbb{C}[M'/M]$. The Borcherds lifting of
 $F$ is then a meromorphic modular form on $\Gamma_0(p)^*$
 whose zeros or poles occur only at cusps or imaginary quadratic
 irrationals corresponding to Heegner points.
 Moreover it has an explicit form of infinite product
 expansion (\cite{B98} Theorem 13.3) which allows us
 to extend (\ref{PRODUCT}) to hold true for all $p\in \mathfrak{S}$.
 (see (\ref{lifting})).
 As a byproduct we can extend Zagier's trace formula to
 all primes $p\in \mathfrak{S}$
 (see Corollary \ref{TRACE}).
 \par
 Corollary \ref{TRACE} gives a description of
 the coefficient of $q^1$ in $f_{d,N}$ as the trace of singular value of
 $j_N^*$.
 In \cite{Zagier}
 \S 7 when the level $N$ is equal to 1 and $D>1$,
 Zagier described the $D$-th coefficient of $f_{d,1}$ in terms of
 the relative
 trace of singular modulus of discriminant $-dD$ from the Hilbert
 class field of $\mathbb{Q}(\sqrt{-dD})$ to its real quadratic
 subfield $\mathbb{Q}(\sqrt{D})$. In \S 4 we will investigate the
 analogue of this in higher level cases and construct Borcherds products
 (see Theorem \ref{trace}).

 \section{Preliminaries}
 \subsection{Modular forms of half-integral weight}
 Here we recall some basic definitions. Let
 $\mathfrak{G}$ be the group consisting of all pairs
 $(A,\psi (\tau))$, where $A=\sm a&b\\c&d \esm\in
 GL_2^+(\mathbb{R})$ and $\psi (\tau)$ is a complex valued
 function holomorphic on $\h$ satisfying $|\psi (\tau)|=(\det
 A)^{-1/4} |c\tau +d|^{1/2}$, with group law defined by
 $(A,\psi_1 (\tau))(B,\psi_2 (\tau))=(AB, \psi_1 (B\tau) \psi_2 (\tau))$.
 If $f:\h \to \mathbb{C}$ and $\xi=(A,\psi (\tau))\in
 \mathfrak{G}$, we put $f|_{[\xi]_{k+1/2}}=f|_\xi=\psi (\tau)^{-2k-1} f(A\tau)$.
 Then from definition it follows that
 $f|_{\xi_1}|_{\xi_2}=f|_{\xi_1 \xi_2}$. There is a monomorphism
 $\Gamma_0(4)\to \mathfrak{G}$ given by $A\mapsto
 A^*:=(A,j(A,\tau))$, where $j(A,\tau)=\left(\frac{c}{d} \right)
  \left(\frac{-1}{d} \right)^{-1/2} (c\tau +d)^{1/2}$
  if $A=\sm a&b\\c&d \esm$.
  A complex valued holomorphic function $f$ on $\h$ is said to be
  a {\em modular form of half-integral weight $k+1/2$ on
  $\Gamma_0(4N)$} if it satisfies $f|_{A^*}=f$ for every $A\in
  \Gamma_0(4N)$, and is holomorphic at the cusps.

 \subsection{Jacobi forms}
 \par
 A {\it (holomorphic) Jacobi form of weight $k$ and index $N$}
 is defined to be a
 holomorphic function $\phi : \mathfrak H \times \mathbb C \to
 \mathbb C$ \, satisfying the two transformation laws
 \begin{align*}
 \phi \left(\frac{a\tau+b}{c\tau+d}, \frac{z}{c\tau+d}\right)
 & = (c\tau+d)^k e^{2\pi i N\frac{cz^2}{c\tau+d}} \phi (\tau,z)
     \q \q (\forall\sm a&b \\c&d \esm \in SL_2(\mathbb Z)), \\
 \phi (\tau, z+\lambda\tau+\mu)
 & = e^{-2\pi iN(\lambda^2 \tau+2\lambda z)} \phi (\tau,z)
     \q \q (\forall\sm \lambda & \mu \esm \in \mathbb Z^2)
 \end{align*}
 and having a Fourier expansion of the form
 \begin{equation}
 \phi (\tau,z)=\sum_{n,r\in \mathbb Z \atop 4Nn-r^2\ge 0} c(n,r) q^n
 \zeta^r \q \q \q (q=e^{2\pi i\tau}, \, \zeta=e^{2\pi iz}),
 \label{Fourier}
 \end{equation}
 where the coefficient $c(n,r)$ depends only on
 $4Nn-r^2$ and on the residue class of $r \, (\text{mod } 2N)$
 (\cite{E-Z} Theorem 2.2).
 In (\ref{Fourier}), if the condition $4Nn-r^2\ge 0$ is deleted,
 we obtain a {\it nearly holomorphic Jacobi form}.
 \par Let $J_{k,N}^!$ be the space of
 nearly holomorphic Jacobi forms of weight $k$ and index $N$.
 Let $J_{*,*}^!$ be the ring of all nearly holomorphic Jacobi
 forms and $J_{ev,*}^!$ its even weight subring. Then $J_{ev,*}^!$
 is the free polynomial algebra over $M_*^!(\Gamma)={\mathbb
 C}[E_4,E_6,\Delta^{-1}]/(E_4^3-E_6^2=1728\Delta)$ on two
 generators $a=\tilde{\phi}_{-2,1}(\tau,z)\in
 {J}_{-2,1}^!$ and $b=\tilde{\phi}_{0,1}(\tau,z)\in
 {J}_{0,1}^!$ (see \cite{E-Z}
 \S 9).
 Let $N\in \{1\} \cup \mathfrak{S}$ and $k=2$.
 There are Jacobi forms $\phi_{D,N}\in J_{2,N}^!$
 uniquely characterized by the requirement that they
 have Fourier
 coefficients $c(n,r)=B(D, 4Nn-r^2)$ which depend only on the
 discriminant $r^2-4Nn$,
 with $B(D,-D)=1$, $B(D,d)=0$ if $d=4Nn-r^2<0,\neq -D$
 and $B(D,0)=-2$ or $0$ according as $D$ is a square or non-square.
 The uniqueness of $\phi_{D,N}$ is obvious since
 the difference
 of any two functions satisfying the definition of $\phi_{D,N}$
 would be an element of $J_{2,N}$, which is of dimension zero
 (see \cite{E-Z} p.118).
 For the existence the
 structure theorem enables us to express $\phi_{D,N}$ as a linear
 combination of $a^i b^{N-i} (i=0,\dots,N)$ over $M_*^!(\Gamma)$.
 \par Through the article we adopt the following notations:
 \begin{itemize}
 \item $\beta$: an element in $\mathbb{Z}/2N\mathbb{Z}$.
 \item $\bar \Gamma'$: the image of $\Gamma'$ in
 $PSL_2(\mathbb{R})$.
 \item $\bar \Gamma'_Q$: $=\{\gamma\in \bar\Gamma' \mid Q\circ \gamma
 = Q\}$.
 \item When we write $\mathcal{Q}_{d,N,\beta}$,
 $\beta$ is always assumed to satisfy $\beta^2\equiv -d$ (mod
 $4N$).
 \item
 $e(x)$: $=\exp (2\pi ix)$.
 \item
 $\zeta_n$: $=e(1/n)=\exp (2\pi i/n)$.
 \item
 $s(D,N)$: the number of prime factors dividing $(D,N)$.
 \item
 $U_m$: the Hecke operator defined by
 $\sum_{n\in\mathbb{Z}}c(n)q^n \mid_{U_m}
 =\sum_{n\in\mathbb{Z}}c(mn)q^n$.
 \item
 {$E_k(\tau)$}: the normalized Eisenstein series of weight $k$,
    equal to $1-(2k/B_{k})\sum_{n>0}\sigma_{k-1}(n)q^n$ where
 $B_k$ is the $k$-th Bernoulli number defined by
  $\sum_{k\in \mathbb{Z}}B_kt^k/k! =t/(e^t-1)$ and
 $\sigma_{k-1}(n)=\sum_{d|n \atop d>0}d^{k-1}$.
 \end{itemize}

\section{Vector valued modular forms}
\begin{Def}
 Let $M$ be an even lattice of signature $(2,n)$
 equipped with a non-degenerate
 quadratic form $q(x)=\frac 12 (x,x)$. Let $M'$ be the
 dual lattice of $M$ and $\mathbb{C} [M'/M]$ be the group ring of
 $M'/M$ with basis $\mathfrak{e}_\gamma$ for each $\gamma\in M'/M$.
 We define a {\em vector valued modular form of weight $k$ and of type
 $\rho_M$} to be a holomorphic function
 $F(\tau)=\sum_{\gamma\in M'/M} h_\gamma(\tau) \e_\gamma$ on the upper half
 plane $\h$ with values in $\mathbb{C} [M'/M]$ such that
 \begin{align*}
    h_\gamma(\tau+1) & = e((\gamma,\gamma)/2) h_\gamma (\tau) \\
    h_\gamma(-1/\tau) & = \frac{\sqrt{i}^{n-2}}{\sqrt{|M'/M|}}
     \sqrt{\tau}^k \sum_{\delta\in M'/M} e(-(\gamma,\delta))
     h_\delta (\tau) \\
    h_{-\gamma}(\tau) & = (-1)^k i^{n-2} h_\gamma(\tau).
\end{align*}
\end{Def}
 For a fixed positive integer $N$, we take $M$ to be the 3
dimensional even lattice of all symmetric matrices $\upsilon= \sm
A&B \\ B&C \esm$ with $NA,B,C\in \mathbb{Z}$ with the norm
$(\upsilon,\upsilon)=-2N\det{(\upsilon)}$. The lattice $M$ splits
as the direct sum of the 2 dimensional hyperbolic unimodular even
lattice
 $II_{1,1}\approx \mathbb{Z} \sm 1/N & 0\\0&0 \esm
 +\mathbb{Z} \sm 0 & 0\\0&1 \esm$ and a lattice
generated by $\sm 0&1\\1&0\esm$ of norm $2N$. Thus the dual
lattice $M'$ is the set of matrices $\sm A&B \\ B&C \esm$ with
$NA,2NB,C\in \mathbb{Z}$, and $M'/M$ can be identified with
$\mathbb{Z}/2N\mathbb{Z}$ by mapping a matrix of $M'$ to the value
of $B\in \frac{1}{2N}\mathbb{Z}/\mathbb{Z}\approx
\mathbb{Z}/2N\mathbb{Z}$. The group $\Gamma_0(N)^*$ acts on the
lattice $M$ as an orthogonal transformation by $\upsilon\mapsto
X\upsilon X^t$ for $\upsilon\in M$ and $X\in \Gamma_0(N)^*$.
Moreover it is not difficult to verify that the special orthogonal
group of $M$ is $\Gamma_0(N)^*$.
 \par
 In the following we will interpret a scalar valued modular form
 in $M_{1/2}^! (N)$ as a vector valued modular form of weight
 $1/2$ and type $\rho_M$. For simplicity we will
 restrict ourselves to the case $N$ is a prime $p$ in the set
 $\mathfrak{S}$.
  Let $f=\sum_{n\in\mathbb{Z}} c(n)
q^n \in M_{1/2}^! (p)$. For $\beta \in \mathbb{Z}/2p\mathbb{Z}$,
we set
$$ h_\beta (\tau) = 2^{s(\beta,p)-1}
\sum_{n\equiv \beta^2 (4p)} c(n) q^{n/4p} $$
 where $s(\beta,p)=1$ if $\beta\equiv 0,p$ (mod $2p$) and
 0 otherwise.
 Then clearly we have
\begin{equation}\label{translation}
h_\beta (\tau+1) = \zeta_{4p}^{\beta^2} h_\beta (\tau).
\end{equation}
For $j$ prime to $4p$, we consider $f_j= f|_{\left(\sm 1&j\\0&4p
\esm, (4p)^{1/4}\right)}|_{W_{4p}}$ where $W_{4p}=\left( \sm
0&-1\\4p&0 \esm,(4p)^{1/4}\sqrt{-i\tau} \right)$. By definition
$f_j=f|_{\left(\sm 4pj&-1\\16p^2&0 \esm, 2\sqrt{-pi\tau}\right)}$.
Let $b$ and $d$ be integers satisfying $jd-4pb=1$ and
$\psi_j=\left(\frac{4p}{j}\right)\sqrt{\left(\frac{-1}{j}\right)}
\zeta_8^{-1}$. Then it holds that $\sm j&b\\4p&d \esm^* \left( \sm
4p&-d \\ 0&4p \esm, \psi_j \right) =\left(\sm 4pj&-1\\16p^2&0
\esm, 2\sqrt{-pi\tau}\right)$. Thus $f_j=f|_{\sm j&b\\4p&d \esm^*
\left( \sm 4p&-d \\ 0&4p \esm, \psi_j \right)}=\psi_j^{-1}
f\left(\tau-\frac{j^{-1}}{4p}\right)$, from which it follows that
\begin{equation}
f|_{\left(\sm 1&j\\0&4p \esm,
(4p)^{1/4}\right)}|_{W_{4p}}=\psi_j^{-1}
f\left(\tau-\frac{j^{-1}}{4p}\right) \label{AA}
\end{equation}
where $j^{-1}$ denotes the inverse of $j$ in
$(\mathbb{Z}/4p\mathbb{Z})^*$. The left hand side of (\ref{AA}) is
\begin{align}
  (4p)^{-1/4} f\left(\frac{\tau+j}{4p}\right)
  |_{W_{4p}} &=(4p)^{-1/4}(\sum_{n\in\mathbb{Z}}
  c(n)\zeta_{4p}^{nj}q^{n/4p})|_{W_{4p}}
  = (4p)^{-1/4}\left( \sum_{\beta (2p)} \zeta_{4p}^{\beta^2j}h_\beta
    \right)|_{W_{4p}} \notag \\
  &=(4p)^{-1/2}(\sqrt{-i\tau})^{-1}
    \left( \sum_{\beta (2p)} \zeta_{4p}^{\beta^2j}h_\beta
    \right)\left(-\frac{1}{4p\tau} \right). \label{BB}
\end{align}
On the other hand, the right hand side of (\ref{AA}) is
\begin{equation}
 \psi_j^{-1}(\sum_{n\in\mathbb{Z}} c(n)\zeta_{4p}^{-nj^{-1}}q^{n})
 =\psi_j^{-1}\left(\sum_{\beta (2p)} \zeta_{4p}^{-\beta^2j^{-1}}h_\beta
    \right)(4p\tau).
\label{CC}
\end{equation}
Replacing $\tau$ by $\tau/4p$ in (\ref{BB}), (\ref{CC}) and
recalling the definition of $\psi_j$ we obtain the following
identity:
\begin{equation}\label{DD}
    \sum_{\beta (2p)} \zeta_{4p}^{\beta^2j}h_\beta(-1/\tau)=\left(
   \frac{4p}{j}\right)\sqrt{\left(\frac{-1}{j}\right)}^{-1}\sqrt{\tau}
   \cdot \sum_{\beta (2p)}
   \zeta_{4p}^{-\beta^2j^{-1}}h_\beta(\tau).
\end{equation}
Let $R$ be a $2p\times 2p$ matrix defined by
$R=\frac{\zeta_8^{-1}}{\sqrt{2p}}\left(
\zeta_{2p}^{-l\beta}\right)_{0\le l,\beta < 2p}.$ To show that
$\sum_{\beta (2p)} h_{\beta} \mathfrak{e}_{\beta} $ is a vector
valued modular form, it is necessary to check that
\begin{equation} \label{inversion}
\sm h_0 \\ \vdots \\ h_\beta \\ \vdots\esm \left( -1/\tau\right)
=\sqrt{\tau} \, R\sm h_0 \\ \vdots \\ h_\beta \\ \vdots\esm
(\tau).
\end{equation}
Since $h_\beta=h_{-\beta}$, the above identity is equivalent to
\begin{equation} \label{EE}
A\sm h_0 \\ \vdots \\ h_\beta \\ \vdots\esm \left( -1/\tau\right)
=\sqrt{\tau} \, AR\sm h_0 \\ \vdots \\ h_\beta \\ \vdots\esm
(\tau)
\end{equation}
for some matrix $A$ with $2p$ columns of which the first $p+1$
ones are linearly independent and span the other ones. We take
 $A=\left(
\zeta_{4p}^{\beta^2 j_l}\right)_{1\le l \le \varphi (4p) \atop
0\le \beta < 2p}$ where $j_l$ is the $l$-th largest element in the
set $\{j \, | \, 1\le j \le 4p \text{ and } (j,4p)=1 \}$. We can
check that the rank of $A$ is $p+1$ unless $p=2$. Moreover from
the Gauss quadratic sum identity (\cite{Lang} Chapter IV, Section
3) it follows that
 $AR=\left( \left(
   \frac{4p}{j_l}\right)\sqrt{\left(\frac{-1}{j_l}\right)}^{-1}
   \zeta_{4p}^{-\beta^2j_l^{-1}} \right)_{1\le l \le \varphi (4p) \atop
 0\le \beta < 2p}$. Now the identity (\ref{EE})
  follows from (\ref{DD}) except for the case $p=2$.

\begin{rmk} When $p=2$, the matrix $A$ is $\sm 1& \zeta_8 & -1 & \zeta_8\\
  1& \zeta_8^3 & -1 & \zeta_8^3\\
 1& \zeta_8^5 & -1 & \zeta_8^5 \\ 1& \zeta_8^7 & -1 & \zeta_8^7\esm$
  and its rank is 2. If we take $j=1,5$ in (\ref{DD}), we obtain the
 following two identities:
 \begin{align}
    h_0\left(-1/\tau\right)+2\zeta_8
    h_1\left(-1/\tau\right)-
    h_2\left(-1/\tau\right)
    &=\sqrt{\tau}(h_0(\tau)+2\zeta_8^{-1} h_1(\tau)-h_2(\tau))
    \label{GG} \\
    h_0\left(-1/\tau\right)-2\zeta_8
    h_1\left(-1/\tau\right)-
    h_2\left(-1/\tau\right)
    &=-\sqrt{\tau}(h_0(\tau)-2\zeta_8^{-1} h_1(\tau)-h_2(\tau)).
    \label{HH}
 \end{align}
 Thus in
 order to determine $h_\beta \left(-1/\tau\right)$ in
 terms of $h_\gamma (\tau)$'s we need one more identity.
 We consider
\begin{align*}
    2(f|_{U_2})|_{(\sm 1&0\\4&1 \esm,\sqrt{4\tau+1} )}
    & =\sum_{j=0}^{1} \sqrt[4]{2}\,
    f|_{(\sm 1&j\\0&2 \esm,\sqrt[4]{2} )(\sm 1&0\\4&1 \esm,\sqrt{4\tau+1}
    )} =\sum_{j=0}^{1} \sqrt[4]{2}\,
    f|_{(\sm 1+4j&j\\8&2 \esm,\sqrt[4]{2}\sqrt{4\tau+1} )}\\
    &=\sqrt[4]{2}\, f|_{\sm 1&0\\8&1 \esm^* (\sm 1&0\\0&2 \esm,\sqrt[4]{2})}
    +
      \sqrt[4]{2}\, f|_{\sm 5&3\\8&5 \esm^* (\sm 1&-1\\0&2 \esm,\sqrt[4]{2})}
    \\
    &=f\left(\frac{\tau}{2} \right)-f\left(\frac{\tau+1}{2}
    \right).
\end{align*}
In terms of  $S=(\sm 0&-1\\1&0 \esm,\sqrt{\tau} )$ and $T=(\sm
1&1\\0&1 \esm,1 )$ we can rewrite the above as
 $2(f|_{U_2})|_{ST^{-4}S^{-1}}=
 f\left(\frac{\tau}{2} \right)-f\left(\frac{\tau+1}{2}\right)$.
 This yields
\begin{equation}\label{FF}
 (f|_{U_2})|_{S}
 =\frac{1}{2}\left(f\left(\frac{\tau}{2} \right)-f\left(\frac{\tau+1}{2}\right)
  \right)|_{ST^4}.
\end{equation}
The left hand side of (\ref{FF}) can be written as
\begin{equation*}
   (h_0(4\tau)+ h_2(4\tau))|_S
    = (h_0(-4/\tau)+h_2(-4/\tau))\sqrt{\tau}^{-1}.
\end{equation*}
Meanwhile the right hand side of (\ref{FF}) turns out to be
\begin{align*}
 2h_1 (4\tau)|_{ST^4}
 &=2(h_1(-4/\tau)\sqrt{\tau}^{-1})|_{T^4} \\
 &=2\left( \frac{\zeta_8^{-1}}{2}\sqrt{\tau/4}
 (h_0(\tau/4)-h_2(\tau/4))\sqrt{\tau}^{-1} \right)|_{T^4}
 \text{ by (\ref{GG}) and (\ref{HH}) } \\
 & =\frac{\zeta_8^{-1}}{2}(h_0(\tau/4)+h_2(\tau/4))
 \text{ by (\ref{translation}). }
\end{align*}
 Replacing $\tau$ by $4\tau$
 we come up with
 \begin{equation}\label{II}
 h_0(-1/\tau)+h_2(-1/\tau) = {\zeta_8^{-1}}\sqrt{\tau}
 (h_0(\tau)+h_2(\tau)).
\end{equation}
Now from (\ref{GG}), (\ref{HH}) and (\ref{II}) we can verify
(\ref{inversion}) in the case $p=2$, too.
\end{rmk}
 Now we will work out Borcherds lifting in $O(2,1)$ case.
 With $M$ as before, we let $G(M)$ denote the {\em Grassmannian} of $M$,
  i.e. the space
 of all two dimensional positive definite subspaces of $M\bigotimes
 \mathbb{R}$ acted on by the orthogonal group $O(2,1)$. In our case the
 Grassmannian $G(M)$ can be identified with the upper half plane
 by mapping $\tau\in \h$ to the two dimensional positive definite
 space spanned by the real and imaginary parts of the norm zero
 vector $\sm \tau^2&\tau\\ \tau&1 \esm$.
 Let $F=\sum_\beta h_\beta \mathfrak{e}_\beta$ be a vector valued
 modular form associated to
 $f=\sum_{n\in\mathbb{Z}}c(n)q^n\in M_{1/2}^! (p)$.
 We write $\sum_{n\in\mathbb{Z}} c_\beta(n/4p)e(n\tau/4p)$
 for the Fourier expansion of $h_\beta(\tau)$.
 Theorem 13.3 of \cite{B98} implies that there is a meromorphic
 function $\Psi_f$ on $G(M)$ that have infinite product expansion
 (explicitly given by Item 5 in Theorem 13.3 of \cite{B98})
 with the following properties:
 \par\noindent
 1. $\Psi_f$ is an automorphic form of weight $c_0(0)/2$ for some
 unitary character of Aut($M,F$).
 \par\noindent
 2. The only zeros or poles of $\Psi_f$ lie on the divisors
 $\lambda^\bot$ for $\lambda\in M$, $\lambda^2 <0$ and are zeros
 of order $\sum_{0< x \in \mathbb{R} \atop x\lambda\in M'}
 c_{x\lambda}(x^2 \lambda^2/2)$.
 \par We can translate the above in terms of usual modular
 forms as follows:
 \begin{thm} \label{language} With the same notation as above,
 we have,
 \par\noindent
 (i) for some $h$,
 $$
 \Psi_f (\tau) = q^{-h} \prod_{n>0} (1-q^n)^{c_\beta(n^2/4p)},
 \text{ (here for each $n>0$, $\beta\equiv n \mbox{ (mod } 2p)$) }
 $$
 is a meromorphic modular form of weight $c_0(0)$ for some
 character of
 $\Gamma_0(p)^*$.
 \par\noindent
 (ii) The possible zeros or poles of $\Psi_f$ occur
 either at cusps or at imaginary quadratic
 irrationals $\alpha\in \h$ satisfying $pA\alpha^2+B\alpha+C=0$ with
 $(A,B,C)=1$.
 Let $\delta=B^2-4pAC<0$. Then the multiplicity of the zero of
 $\Psi_f$ at $\alpha$ is given by $\sum_{n>0} c_\beta(\delta n^2/4p)$
 (here
 for each $n>0$,
 $\beta$ is chosen such that
 $\beta^2\equiv\delta n^2  \mbox{ (mod $4p$)})$.
 \end{thm}
  Recall the nearly holomorphic modular form
 $f_{d,p}$ in $M_{1/2}^!(p)$ with Fourier expansion of the form
 $q^{-d}+\sum_{D\ge 1} A(D,d)q^d$.
 We write
  $2f_{d,p}=\sum_{n\in\mathbb{Z}} c(n) q^n$.
  By Theorem \ref{language} (i),
  \begin{equation} \label{first}
  \Psi_{2f_{d,p}}=q^{-h} \prod_{u=1}^\infty
 (1-q^u)^{A^*(u^2,d)} \q\q\q \mbox{(for some } h)
  \end{equation}
  is a meromorphic modular form of weight 0 for
  some character of $\Gamma_0(p)^*$.
 Choose an integer $\beta \, (\mbox{mod }
 2p)$ with $\beta^2 \equiv -d \, (\mbox{mod } 4p)$.
  From Theorem \ref{language} (ii) it follows that the zeros of
  $\Psi_{2f_{d,p}}$ occur at $\alpha_Q$
  for each $Q\in \mathcal{Q}_{d,p,\beta}$ with multiplicity
  $$ c_\beta(-d/4p)=2^{s(\beta,p)-1} c(-d) =\left\{%
\begin{array}{ll}
    1 & \Leftrightarrow p\nmid \beta
    \Leftrightarrow \beta\not\equiv -\beta \mbox{ (mod $2p$)}
    \Leftrightarrow Q \circ W_p \not\in  \mathcal{Q}_{d,p,\beta} \\
    2 &  \Leftrightarrow p\mid \beta
    \Leftrightarrow \beta\equiv -\beta \mbox{ (mod $2p$)}
    \Leftrightarrow Q \circ W_p \in  \mathcal{Q}_{d,p,\beta}.
\end{array}%
\right. $$
 On the other hand, the multiplicity of the zero of
 $j_p^*(\tau)-j_p^*(\alpha_Q)$ at $\alpha_Q$ is
 $$ |\bar\Gamma_0(p)^*_Q|=
\left\{%
\begin{array}{ll}
    |\bar\Gamma_0(p)_Q|, & \hbox{ if $Q \circ W_p\neq Q$
     (mod $\Gamma_0(p)$); } \\
    2|\bar\Gamma_0(p)_Q|, & \hbox{ if $Q \circ W_p=Q$ (mod $\Gamma_0(p))$. } \\
\end{array}%
\right.
 $$
 From this observation we come up with
 \begin{equation} \label{second}
 \Psi_{2f_{d,p}} =\prod_{Q\in \mathcal{Q}_{d,p,\beta}/\Gamma_0(p)}
  (j_p^*(\tau)-j_p^*(\alpha_Q))^{\frac{1}{|\bar{\Gamma}_0(p)_Q|}}
 \end{equation}
 Comparing (\ref{first}) with (\ref{second}) we get
 \begin{equation} \label{lifting}
 q^{-H_{p,\beta}(d)}\prod_{u=1}^\infty
 (1-q^u)^{A^*(u^2,d)}
 =\prod_{Q\in \mathcal{Q}_{d,p,\beta}/\Gamma_0(p)}
  (j_p^*(\tau)-j_p^*(\alpha_Q))^{\frac{1}{|\bar\Gamma_0(p)_Q|}}
 \end{equation}
 where $H_{p,\beta}(d)=\sum_{Q\in \mathcal{Q}_{d,p,\beta}/\Gamma_0(p)}
 \frac{1}{|\bar\Gamma_0(p)_Q|}$.
 Note that if $d$ is not divisible as a discriminant by
 $p^2$, then $|\mathcal{Q}_{d,p,\beta}/\Gamma_0(p)|
 =|\mathcal{Q}_{d}/\Gamma|$ and $|\bar\Gamma_0(p)_Q|=|\bar\Gamma_Q|$.
 Therefore in this case $H_{p,\beta}(d)$ is none other than
 $H(d)$.
 \begin{Exm}
 Let $p=2, d=16$ and $\beta\equiv 0 \mbox{ (mod 4)}$.
 A set of representatives for
 $ \mathcal{Q}_{16,2,0}/\Gamma_0(2)$ is given by
 $Q_1=[4,-4,2]$, $Q_2=[2,0,2]$ and $Q_3=[4,0,1]$. We observe that
 $ |\bar\Gamma_0(2)_{Q_i}|=2$ if $i=1$ and 1 otherwise.
 Up to modulo $\Gamma_0(2)$, the involution $W_2$ fixes $Q_1$ and
 interchanges $Q_2$ with $Q_3$. Thus, by (\ref{second}),
 $\Psi_{2f_{16,2}}=(j_2^*(\tau)-j_2^*(\alpha_{Q_1}))^{1/2}
 (j_2^*(\tau)-j_2^*(\alpha_{Q_2}))^2$.
 \end{Exm}
 Let $f=\sum_{D\in \mathbf{Z}}A(D) q^D \in
 M^!_{1/2}(p)$. As before ,
 we set
 $ h_\beta (\tau) = 2^{s(\beta,p)-1}
 \underset{D\equiv \beta^2 (4p)}{\sum} A(D) q^{D/4p}$ for each $\beta\in
\mathbb{Z}/2p\mathbb{Z}$. By (\ref{inversion})
 $$
 h_\beta (-1/\tau)=\frac{\sqrt{i}^{-1}\sqrt{\tau}}{\sqrt{2p}}
 \sum_{j \, (2p)} e(-\beta j/2p)h_j(\tau).
 $$
 Let $\phi\in J^!_{2,p}$. The Fourier coefficient $c(n,r)$ of
 $q^n\zeta^r$ in $\phi$ depends only on the discriminant $r^2-4pn$
 and therefore we can write $c(n,r)=B(4pn-r^2)$.
 For $\mu\in \mathbb{Z}/2p\mathbb{Z}$, we set
 $\tilde{h}_\mu=\sum_{d\equiv -\mu^2 (4p)} B(d) q^{d/4p}.$
 The formula (9) in \cite{E-Z} \S 5 says
 $$
 \tilde{h}_\mu (-1/\tau)=\frac{\sqrt{i}\sqrt[3]{\tau}}{\sqrt{2p}}
 \sum_{k \, (2p)} e(\mu k/2p)\tilde{h}_k(\tau).
 $$
 Let $g=\sum_{d\in\mathbb{Z}} B(d) q^d$ with $B(d)$ defined to be
 zero unless $-d$ is a square modulo $4p$. From definition
 it holds that
 $f=\sum_{\beta \, (2p)} h_\beta (4p\tau)$ and
 $g=\sum_{\mu \, ( 2p)} 2^{s(\mu,p)-1}\tilde{h}_\mu (4p\tau)$.
 Moreover $(fg)|_{U_{4p}}=\sum_{i \, ( 2p)}
 (h_i\tilde{h}_i)(\tau)$.
 Thus
 $(fg)|_{U_{4p}}(\tau+1)=(fg)|_{U_{4p}}(\tau)$ and
 \begin{align*}
    (fg)|_{U_{4p}}(-1/\tau) & =\sum_{i \, (2p)}
 (h_i\tilde{h}_i) (-1/\tau)
  = \frac{\tau^2}{2p} \sum_{i,j,k \, (2p)}
  e(-ij/2p)e(ik/2p)(h_j\tilde{h}_k)(\tau) \\
  &= \frac{\tau^2}{2p} \sum_{j,k \, (2p)}
   \sum_{i \, (2p)} e(i(-j+k)/2p)(h_j\tilde{h}_k)(\tau) \\
  &= \tau^2 \sum_{j \, (2p)} (h_j\tilde{h}_j)(\tau)
   \text{ \, since } \sum_{i \, (2p)} e(i(-j+k)/2p)=
  \begin{cases}
    2p, & \hbox{if $j\equiv k$ (mod $2p$)} \\
    0, & \hbox{otherwise}
  \end{cases} \\
  &=\tau^2 (fg)|_{U_{4p}}.
\end{align*}
 $(fg)|_{U_{4p}}$ is then invariant under the action of
 $SL_2(\mathbb{Z})$ with a pole only at $\infty$ and of weight 2. Thus
 $(fg)|_{U_{4p}}$ can be written as the derivative of some
 polynomial in $j$. In particular
 the constant term of the $q$-expansion of $fg$ vanishes.
 Applying this to $f=f_{d,p}$ and $\phi=\phi_{D,p}$ gives rise to
 the following lemma:

 \begin{lem} \label{duality}
 Let $d=4pn-r^2$ for some integers $n$ and $r$. Let
 $A(D,d)$ (resp. $B(D,d)$)be the coefficient of $q^D$
 (resp. $q^n \zeta^r$) in
 $f_{d,p}\in M^!_{1/2}(p)$ (resp. $\phi_{D,p}\in J^!_{2,p}$). Then
 $A(D,d)=-B(D,d)$.
 \end{lem}
 \begin{cor} \label{TRACE}
 Let $p$ be a prime for which $\Gamma_0(p)^*$ is of genus zero.
 For each natural number $d$ which is congruent to a square modulo
 $4p$, let
 $\textbf{\emph{t}}^{(p)}(d)=\sum_{Q\in\mathcal{Q}_{d,p,\beta}/\Gamma_0(p)}
  \frac{1}{|\bar\Gamma_0(p)_Q|} j_p^* (\alpha_Q)$. We also put
 $\textbf{\emph{t}}^{(p)}(-1)=-1$, $\textbf{\emph{t}}^{(p)}(0)=2$ and
 $\textbf{\emph{t}}^{(p)}(d)=0$
 for $d<-1$. Then the series
 $\sum_{n,r} \textbf{\emph{t}}^{(p)}(4pn-r^2) q^n \zeta^r$ becomes a
 nearly holomorphic Jacobi form of weight 2 and index $p$.
 \end{cor}
 \begin{proof}
 In (\ref{lifting}) if we compare the $q$-expansions on both
 sides, then we obtain $\textbf{t}^{(p)}(d)=A^*(1,d)=A(1,d)$ which
 is equal to $-B(1,d)$ by Lemma \ref{duality}. Now the series
 $\sum_{n,r} \textbf{t}^{(p)}(4pn-r^2) q^n \zeta^r$ coincides with
 $-\phi_{1,p}$.
 \end{proof}

\section{Twisted traces of singular values}

 Let $D$ and $-d$ be a positive and a negative discriminant,
 respectively. For simplicity we suppose them to be
 fundamental and coprime.
 Let $K=\mathbb Q(\sqrt{-Dd})$.
 Let $p_1,\dots, p_r$ be odd primes dividing $Dd$ and
 put $p_i^*=(-1)^{(p_i-1)/2} p_i$ for $i=1,\dots,r$.
 By class field theory
 the genus field $M$ of $K$ is expressed as $K(\sqrt{p_1^*},\dots,
 \sqrt{p_r^*})$ and is contained in the Hilbert class field $H$ of
 $K$ (see \cite{Cox} Theorem 6.1).
 We have obvious inclusion $N=K(\sqrt{D})\subset
 M\subset H$.
 The genus character $\chi=\chi_{D,-d}$ assigns to any quadratic
 form $Q$ of discriminant $-dD$ a value $\pm 1$ defined by
 $\chi (Q)=\left(\frac{D}{q} \right)=\left(\frac{-d}{q} \right)$
 where $q$ is any prime represented by $Q$ and not dividing $Dd$.
 This is independent of the choice of $q$.
 We further assume $D$ and $-d$ to be congruent to a square modulo $4p$.
 We define the ``twisted trace" $\textbf{t}^{(p)}(D,d)$ by
 $$\textbf{t}^{(p)}(D,d)=\sum_{Q\in \mathcal{Q}_{dD,p,\beta}/\Gamma_0(p)} \chi (Q)
 j_p^*(\alpha_Q).$$
 \begin{prop} \label{twisted}
 {\emph{(i)}} The definition of \,
 $\textbf{\emph{t}}^{(p)}(D,d)$ is independent of the choice of
 $\beta$.
 \par\noindent
 {\emph{(ii)}} Let $p$ be a prime such that $p-1$ divides $24$. Then
 $\frac{1}{\sqrt{D}} \textbf{\emph{t}}^{(p)}(D,d)$ is a rational integer.
 \end{prop}
 \begin{proof}
 (i) Since $W_p$
 induces a bijection between
 ${ \mathcal{Q}_{dD,p,\beta}/\Gamma_0(p)}$ and
 ${ \mathcal{Q}_{dD,p,-\beta}/\Gamma_0(p)}$,
 it is enough to show that
  $\chi (Q\circ W_p)= \chi (Q)$ and $j_p^*(\alpha_{Q\circ W_p
  })=j_p^*(\alpha_Q)$ for each $Q\in \mathcal{Q}_{dD,p,\beta}$.
  The latter is clear since $j_p^*$ is the Hauptmodul for
  $\Gamma_0(p)^*$. The invariance of $\chi$ under $W_p$ follows
  from Proposition 1 in p.508 of \cite{G-K-Z}.
  \par\noindent
  (ii) It is well-known that $j(\alpha_Q)$ is an algebraic
  integer. By \cite{Chen-Yui} Remark 1.5.3 and Appendix 1 there is a
  modular relation between $j$ and $j_p^*$ which allows $j_p^*$ to
  be integral over $\mathbb{Z}[j]$. Therefore $j_p^*(\alpha_Q)$ is
  an algebraic integer, too. Let $j_{0,p}$ be the Hauptmodul for
  $\Gamma_0(p)$. Then $j_p^*$ can be written as $j_{0,p}(\tau) +
  j_{0,p}(W_p \tau)$ up to some constant ($\in \mathbb{Z}$).
  By \cite{Chen-Yui} Theorem 3.7.5,
  $j_{0,p}(\alpha_Q)$ (resp. $j_{0,p}(W_p \alpha_Q)$)
   generates Hilbert class
  field $H$ of $K$ so that $j_p^*(\alpha_Q)$ belongs to $H$.
  The set of forms $Q$ with $\chi
  (Q)=+1$ is identified with ${\rm Gal}(H/N)$.
  Therefore the sum $\sum_{\chi (Q) =+1} j_p^*(\alpha_Q)$ is the trace
  of the algebraic integer $j_p^*(\alpha_Q)$ from $H$ to $N$. Since
  this number is invariant under the complex conjugation, it
  belongs to $\mathbb{Q} (\sqrt{D})$. If we write this number as
  $\lambda=\frac{1}{2} (a+b\sqrt{D})$ with $a,b\in \mathbb{Z}$,
  then the sum over all $j_p^*(\alpha_Q)$ is equal to
  ${\rm Tr}_{H/K} (j_p^*(\alpha_Q))={\rm Tr}_{N/K}\circ {\rm Tr}_{H/N}
   (j_p^*(\alpha_Q))={\rm Tr}_{N/K}(\lambda)=a$. Hence the sum
   $\textbf{t}^{(p)}(D,d)$ is equal to $-a+2\lambda=b\sqrt{D}$.
 \end{proof}
 Here are some numerical examples with $p=2$ and $p=3$:
 \par\noindent
 $
 \frac{1}{\sqrt{17}}\textbf{t}^{(2)}(17,4)=\frac{1}{\sqrt{17}}(
 j_2^*\left( \alpha_{ [18,2,1]}\right)+
 j_2^*\left( \alpha_{ [2,2,9]}\right)-
 j_2^*\left( \alpha_{ [6,-2,3]}\right)-
 j_2^*\left( \alpha_{ [6,2,3]}\right))= -204800$,
 \par\noindent
 $
 \frac{1}{\sqrt{8}}\textbf{t}^{(2)}(8,7)=\frac{1}{\sqrt{8}}(
 j_2^*\left( \alpha_{ [14,0,1]}\right)+
 j_2^*\left( \alpha_{ [2,0,7]}\right)-
 j_2^*\left( \alpha_{ [6,4,3]}\right)-
 j_2^*\left( \alpha_{ [10,8,3]}\right))= 90112$,
 \par\noindent and
 \par\noindent
 $ \frac{1}{\sqrt{13}} \textbf{t}^{(3)}(13,3)=\frac{1}{\sqrt{13}}
 (j_3^*\left( \alpha_{[12,3,1] }\right)+
 j_3^*\left( \alpha_{ [3,3,4]}\right)-
 j_3^*\left( \alpha_{[6,3,2] }\right)-
 j_3^*\left( \alpha_{[15,9,2] }\right)) = -378,$
 \par\noindent
 $
 \frac{1}{\sqrt{13}}\textbf{t}^{(3)}(13,8)=\frac{1}{\sqrt{13}}(
 j_3^*\left( \alpha_{ [27,2,1]}\right)-
 j_3^*\left( \alpha_{ [21,8,2]}\right)+
 j_3^*\left( \alpha_{ [3,2,9]}\right)+
 j_3^*\left( \alpha_{ [9,2,3]}\right)-
 j_3^*\left( \alpha_{ [6,-4,5]}\right)-
 j_3^*\left( \alpha_{ [15,14,5]}\right))$ $=-11968 $,
 \par\noindent
 $
 \frac{1}{\sqrt{21}}\textbf{t}^{(3)}(21,8)=\frac{1}{\sqrt{21}}(
 j_3^*\left( \alpha_{ [42,0,1]}\right)-
 j_3^*\left( \alpha_{ [21,0,2]}\right)+
 j_3^*\left( \alpha_{ [3,0,14]}\right)-
 j_3^*\left( \alpha_{ [6,0,7]}\right))= 342144$.
 \par\noindent

 We want to express these traces as the coefficients of $f_{d,p}$.
 As explained in the  Appendix of \cite{Kim}, $f_{d,p}$ can be found by
 making use of
 ``Rankin-Cohen bracket" $[ \q , \q ]_n$.
 Let $f(\tau)$ and $g(\tau)$ denote two modular forms of weight
 $k$ and $l$ on some group $\Gamma'\subset \Gamma$. We denote by
 $D$ the differential operator $\frac{1}{2\pi i} \frac{d}{d\tau}=
 q\frac{d}{dq}$ and use $f',f'',\dots,f^{(n)}$ instead of
 $Df,D^2f,\dots,D^n f$. The $n$-th Rankin-Cohen bracket of $f$
 and $g$ is defined by the formula
 $$ [f,g]_n (\tau)=\sum_{r+s=n} (-1)^r \mat n+k-1\\s\emat
 \left(\begin{matrix} n+l-1\\r\end{matrix}\right) f^{(r)}(\tau) g^{(s)}(\tau).$$
 By \cite{Cohen} \S 7 (or \cite{Zagier94} \S 1) $[f,g]_n (\tau) $
 is a modular form of weight $k+l+2n$ on $\Gamma'$.
 Let us consider the construction of $f_{d,p}$ in the cases $p=2$ and $p=3$:
 \par\noindent
 $p=2$:
 Let $\theta=\sum_{n\in \mathbb Z} q^{n^2}$,
 $u=([\theta,E_{10}(8\tau)]_1/\Delta(8\tau)+2112\theta)/(-40)$ and
 $v=([\theta,E_{8}(8\tau)]_2/\Delta(8\tau)-11520\theta)/72$.
 Then it follows that
 \par\noindent
 \par\noindent
 $f_{4,2}=(v-u)/12=q^{-4}-52q+272q^4+2600q^8-8244q^9
 +15300q^{12}+71552q^{16}-204800q^{17}+282880q^{20}+\cdots,$
 \par\noindent
 $f_{7,2}=(4u-v)/3=q^{-7}-23q-2048q^4+45056q^8+252q^9
 -516096q^{12}+4145152q^{16}-1771q^{17}-26378240q^{20}+\cdots$

 \par\noindent
 $p=3$: Let
 $u=([\theta,E_{10}(12\tau)]_1/\Delta(12\tau)+1584\theta)/(-20)$,
 $v=([\theta,E_{8}(12\tau)]_2/\Delta(12\tau)-25920\theta)/72$ and
 $w=([\theta,E_{6}(12\tau)]_3/\Delta(12\tau)+272160\theta)/(-112)$.
 Then we obtain that
 \par\noindent
 $f_{3,3}=(4u-5v+w)/360=q^{-3}-14q+40q^4-78q^9
 +168q^{12}-378q^{13}+688q^{16}-897 q^{21}\cdots,$
 \par\noindent
 $f_{8,3}=(-9u+10v-w)/60=q^{-8}-34q-188q^4+2430q^9
 +8262q^{12}-11968q^{13}-34936q^{16}+171072 q^{21}\cdots,$
 \par\noindent
 $f_{11,3}=(36u-13v+w)/24=q^{-11}+22q-552q^4-11178q^9
 +48600q^{12}+76175q^{13}-269744q^{16}
 -1782891 q^{21}\cdots.$

 Comparing these coefficients with the examples of twisted traces
 $\textbf{t}^{(p)}(D,d)$ we are led to the following result which
 can be proved in a way analogous to the proof of \cite{Zagier}
 Theorem 6:
 \begin{equation} \label{Twisted_Trace}
 \frac{1}{\sqrt{D}}\textbf{t}^{(p)}(D,d)=A^*(D,d)=-B^*(D,d)
 \end{equation}
 with $A^*(D,d)=2^{s(D,p)}A(D,d)$ and $B^*(D,d)=2^{s(D,p)}B(D,d)$.

 Now we will derive certain product identities from
 (\ref{Twisted_Trace}).
 \begin{lem}
 For $i\ge 0$ and $m$ coprime to $p$,
 $$
 \phi_{p^{2i}m^2D, 1}|_{V_p}=p\phi_{p^{2i+2}m^2D,p}+\phi_{p^{2i}m^2D,p}.$$
 Here $V_p$ is the Hecke operator on Jacobi forms defined by the
 formula (2) in \cite{E-Z}.
 \label{Hecke1}
 \end{lem}
 \begin{proof}
 According to \cite{E-Z} Theorem 4.1, the operator $V_p$ maps
 $J_{2,1}^!$ to $J_{2,p}^!$.
 From the formula (7) in \cite{E-Z}, p.43, we find that
 $$
 \text{the coefficient of $q^n\zeta^r$ in
 }\phi_{p^{2i}m^2D,1}|_{V_p}=
 \begin{cases} p, & \text{ if } 4pn-r^2=-p^{2i+2}m^2D \\
                1, & \text{ if } 4pn-r^2=-p^{2i}m^2D \\
                0, & \text{ if } 4pn-r^2<0,
                    \neq -p^{2i+2}m^2D,-p^{2i}m^2D.\end{cases}$$
 From these observations and the uniqueness of $\phi_{D,p}$, the lemma
 immediately follows.
 \end{proof}
 Let $t=j_p^*$.
 For each positive integer $m$, we define
 $$\textbf{t}_m^{(p)}(D, d)=\sum_{Q\in
 {\mathcal Q}_{dD,p,\beta}/\Gamma_0(p)}
  \chi (Q) t_m(\alpha_Q)$$
 where $t_m$ is a unique polynomial of $t$ such that $t_m\equiv
 q^{-m}$ (mod $q\mathbb{C}[[q]]$).
 We can use Hecke operators in integral and half-integral weight
 to get the following lemma whose proof is similar to that of
 formulas in p.13 of \cite{Zagier}.
 \begin{lem}
 Let $l$ be a positive integer coprime to $p$ and $d=4pn-r^2$. Then
 \par\noindent
 (i) $\frac{1}{\sqrt{D}}\textbf{\emph{t}}_l^{(p)}(D, d)=
 -\text{coefficient of $q^n\zeta^r$ in }2^{s(D,p)}\phi_{D,p}|_{T_l}$,
 \par\noindent
 (ii) $\phi_{D,p}|_{T_l}=\sum_{\nu | l} \left(\frac{D}{l/\nu} \right)
 \nu \phi_{\nu^2D,p}$.
 \label{Hecke2}
 \end{lem}
 We claim that for $d=4pn-r^2$,
 \begin{equation}
  \frac{1}{\sqrt{D}}\textbf{t}_m^{(p)}(D, d)
  =-\text{coefficient of $q^n\zeta^r$ in }
  \sum_{u|m} 2^{s(u^2D,p)}\left(\frac{D}{m/u} \right)u\phi_{u^2D,p}.
  \label{claim}
 \end{equation}
 Let $J=j-744$ and $J_m$ be the unique polynomial
 such that $J_m\equiv
 q^{-m}$ (mod $q\mathbb{C}[[q]]$).
 By ``expansion formula" (\cite{Koike} p.14 or \cite{C-N} p.319)
 we come up with $t_{lp^{k+1}}(\tau)=J_{lp^{k}}(p\tau)+J_{lp^{k}}(\tau)
  -t_{lp^{k}}(\tau)$ (for $k\ge 0$) and therefore
 \begin{align}
 &\sum_{Q\in {\mathcal Q}_{dD,p,\beta}/\Gamma_0(p)} \chi (Q)
 t_{lp^{k+1}}(\alpha_Q) \notag \\
 =&\sum_{Q\in {\mathcal Q}_{dD,p,\beta}/\Gamma_0(p)}
  \chi (Q) (J_{lp^{k}}(p\tau)+J_{lp^{k}}(\tau))|_{\tau
  =\alpha_Q}
 -\sum_{Q\in {\mathcal Q}_{dD,p,\beta}/\Gamma_0(p)} \chi (Q)
 t_{lp^{k}}(\alpha_Q).
 \label{Hecke}
 \end{align}
 The map which sends $[a,b,c]\in {\mathcal Q}_{dD,p,\beta}$ to
 $[a/p,b,cp]\in
 {\mathcal Q}_{dD}$ induces a bijection
 between ${\mathcal Q}_{dD,p,\beta}/\Gamma_0(p)$
 and ${\mathcal Q}_{dD}/\Gamma$. And the natural map from
 ${\mathcal Q}_{dD,p,\beta}/\Gamma_0(p)$ to ${\mathcal Q}_{dD}/\Gamma$
 also gives
 a bijection. Theses maps preserve the values of $\chi$ and
 (\ref{Hecke}) is rewritten as
 $$
 \sum_{Q\in {\mathcal Q}_{dD,p,\beta}/\Gamma_0(p)} \chi (Q)
 t_{lp^{k+1}}(\alpha_Q)
 =2\sum_{Q\in {\mathcal Q}_{dD}/\Gamma}
  \chi (Q) J_{lp^{k}}(\alpha_Q)
 -\sum_{Q\in {\mathcal Q}_{dD,p,\beta}/\Gamma_0(p)} \chi (Q)
 t_{lp^{k}}(\alpha_Q),
 $$ which yields
 \begin{equation*}
 \textbf{t}_{lp^{k+1}}^{(p)}(D,d)=2\textbf{t}_{lp^{k}}^{(1)}(D,d)
 -\textbf{t}_{lp^{k}}^{(p)}(D,d)
 \text{ \hspace{0.5cm} for }k\ge 0.
 \end{equation*}

 In (\ref{claim}) we write $m=l p^k$ with $(l,p)=1$.
 We use induction on $k$ to prove the claim.
 If $k=0$, the claim (\ref{claim}) follows from Lemma \ref{Hecke2}.
 Now we assume the claim for $k$.
 \begin{align*}
 & \frac{1}{\sqrt{D}} \textbf{t}_{lp^{k+1}}^{(p)}(D,d)
  = \frac{1}{\sqrt{D}}(2\textbf{t}_{lp^{k}}^{(1)}(D,d)
  -\textbf{t}_{lp^{k}}^{(p)}(D,d)) \\
 &=-\text{coefficient of $q^n\zeta^r$ in }
   2(\phi_{D,1}|_{T_{lp^k}})|_{V_p}
  -\frac{1}{\sqrt{D}}\textbf{t}_{lp^{k}}^{(p)}(D,d)
  \text{  by \cite{Zagier} formula (25), Theorem 5-(iii)
 }
 \\
 &=-\text{coefficient of $q^n\zeta^r$ in }
  2(\sum_{i=0}^k\sum_{\nu|l}\nu p^i\left( \frac{D}{(l/\nu)p^{k-i}}\right)
  \phi_{\nu^2p^{2i}D,1})|_{V_p}
  -\frac{1}{\sqrt{D}}\textbf{t}_{lp^{k}}^{(p)}(D,d)\\
 & \text{ \hspace{1cm} by \cite{Zagier} p.13 and Theorem 5-(iii)
  }
 \\
 &=-\text{coefficient of $q^n\zeta^r$ in } [
  2\sum_{i=0}^k\sum_{\nu|l}\nu p^i\left( \frac{D}{(l/\nu)p^{k-i}}\right)
  (p\phi_{\nu^2p^{2i+2}D,p}+\phi_{\nu^2p^{2i}D,p})
  \\
 & \hspace{1cm} \q -(\sum_{i=1}^k\sum_{\nu|l}
  2\nu p^i\left( \frac{D}{(l/\nu)p^{k-i}}\right)\phi_{\nu^2p^{2i}D,p}
  +\sum_{\nu|l}2^{s(D,p)}\left( \frac{D}{(l/\nu)p^{k}}\right)
  \nu \phi_{\nu^2D,p}) ] \\
 &  \text{ \hspace{1cm} by Lemma \ref{Hecke1} and induction hypothesis} \\
 &=-\text{coefficient of $q^n\zeta^r$ in }
  [ \sum_{i=1}^{k+1}\sum_{\nu|l}
  2\nu p^i\left( \frac{D}{(l/\nu)p^{k+1-i}}\right)
 \phi_{\nu^2p^{2i}D,p}
 +\sum_{\nu|l}2^{s(D,p)}\left( \frac{D}{(l/\nu)p^{k+1}}\right)
 \nu\phi_{\nu^2D,p} ]  \\
  & \text{ \hspace{1cm} since }
  \left(\frac{D}{p} \right)=
   0 \text{ or } 1, \text{ and }
  2^{s(D,p)}=1 \text{ or } 2.
 \end{align*}

 Let
 $z\in {\mathfrak H}$.
 Note that $\frac 1m t_m(z)$ can be viewed as the coefficient of $q^m$-term
 in $-\log q - \log (t(\tau)-t(z))$ (see \cite{Norton}). Thus
 $\log q^{-1}-\sum_{m > 0} \frac 1m t_m(z) q^m = \log (t(\tau)-t(z))$.
 Taking exponential on both sides, we get
 \begin{equation}
 q^{-1} \exp (-\sum_{m> 0} \frac 1m t_m(z) q^m)=t(\tau)-t(z).
 \label{product1}
 \end{equation}
 By the claim (\ref{claim}), we obtain
 \begin{equation}
 \frac{1}{\sqrt{D}} \textbf{t}^{(p)}_m(D,d)=-\sum_{u|m}u \left(\frac{D}{m/u} \right)
  B^*(u^2D,d).
 \label{Hecke3}
 \end{equation}
  From (\ref{product1}) and (\ref{Hecke3}) it follows that
 \begin{align*}
  & \prod_{Q\in \mathcal{Q}_{dD,p,\beta}/\Gamma_0(p)}
  (t(\tau)-t(\alpha_Q))^{\chi(Q)}
  =\exp (-\sum_{m=1}^\infty \textbf{t}^{(p)}_m(D,d) q^m/m) \\
  &=\exp (\sqrt{D}\sum_{m=1}^\infty \sum_{u|m}
  u \left(\frac{D}{m/u} \right) B^*(u^2D,d) q^m/m)
  =\exp (\sqrt{D}\sum_{m=1}^\infty \sum_{u=1}^\infty
   u\left(\frac{D}{m} \right)B^*(u^2D,d) \frac{q^{mu}}{mu}) \\
  &=\exp (\sum_{u=1}^\infty (-B^*(u^2D,d))
   (-\sum_{0<n<D}\left(\frac{D}{n} \right)
     \sum_{m=1}^\infty (\zeta_D^n q^u)^m/m))
   \\
 & \text{ \hspace{1cm} by the Gauss sum identity }
 \sqrt{D}\left(\frac{D}{m} \right)
 =\sum_{0<n<D}\left(\frac{D}{n} \right) \zeta_D^{nm}
 \\
 &=\exp (\sum_{u=1}^\infty A^*(u^2D,d)
   \sum_{0<n<D}\left(\frac{D}{n} \right)
     \log (1-\zeta_D^n q^u)).
 \end{align*}
 Now the following theorem is immediate.
 \begin{thm}\label{trace}
 Let $D>1$ and $-d<0$ be coprime fundamental discriminants and
 congruent to a square modulo $4p$.
 We define
 $
 \mathcal{H}_{D,d,p}(X)=\prod_{Q\in \mathcal{Q}_{dD,p,\beta}/\Gamma_0(p)}
  (X-j_p^*(\alpha_Q))^{\chi(Q)}.
 $
 Then
 $$
 \mathcal{H}_{D,d,p}(j_p^*(\tau))=\prod_{u=1}^{\infty} P_D(q^u)^{A^*(u^2D,d)}
 $$
 where $A^*(D,d)=2^{s(D,p)}A(D,d)$ and
 $P_D(t)=\prod_{0<n<D} (1-\zeta_D^n t)^{\left(\frac{D}{n}\right)}$.
 \end{thm}

 \begin{Exm}
 Let us consider the case $p=3$ and $D=13, d=3$.
 For $\beta\equiv 1 \mbox{ (mod 6)}$, a set of representatives for
 $\mathcal{Q}_{39,3,1}/\Gamma_0(3)$ is given by
 $ \{ [12,3,1],[3,3,4],[6,3,2],[15,9,2]\}. $
 Meanwhile the $q$-expansion of $f_{3,3}$ is
 \begin{equation*}
 q^{-3}-14q+40q^4-78q^9+168q^{12}-378q^{13}+688q^{16}
 +\cdots+133056q^{52}+\cdots+(-30650256)q^{117}+\cdots.
 \end{equation*}
 Theorem \ref{trace} yields the following product
 formula:
 \begin{equation*}
  \frac{\left(j_3^*(\tau)-j_3^*\left( \alpha_{[12,3,1] }\right)\right)
   \left(j_3^*(\tau)-j_3^*\left( \alpha_{[3,3,4] }\right)\right)}
   {\left(j_3^*(\tau)-j_3^*\left( \alpha_{[6,3,2] }\right)\right)
    \left(j_3^*(\tau)-j_3^*\left( \alpha_{[15,9,2] }\right)\right) }
  =P_{13}(q)^{-378}P_{13}(q^2)^{133056}P_{13}(q^3)^{2(-30650256)}\cdots.
 \end{equation*}
 \end{Exm}

 \begin{center}
{\bf Acknowledgment}
\end{center}
 I would like to thank Professor J. H. Bruinier for encouraging me to
 do this subject and kindly answering my questions.
 I also want to express my gratitude to the people at
 Mathematisches Institut der Universit\"{a}t Heidelberg  for their
 support.


\begin{thebibliography}{99}
 \bibitem{B92}
       R. E. Borcherds,
       {\em Monstrous moonshine and monstrous Lie superalgebras},
       Invent. math. {\bf 109} (1992), 405-444.

 \bibitem{B95}
  R. Borcherds,
  {\em Automorphic forms on $O_{s+2,2}(\mathbb R)$ and infinite
  products},
  Invent. Math. {\bf 120} (1995),
  161-213.

  \bibitem{B98}
  R. Borcherds,
  {\em Automorphic forms with singularities on Grassmanians},
  Invent. Math. {\bf 132} (1998),
  491-562.

  \bibitem{B99}
  R. Borcherds,
  {\em The Gross-Kohnen-Zagier theorem in higher dimensions},
  Duke Math. J. {\bf 97} (1999),
  219-233.

  \bibitem{Bruinier}
  J. H. Bruinier,
  {\em Borcherds products on $O(2,l)$ and Chern classes of Heegner
  divisors},
  Springer Lecture Notes in Mathematics 1780, Springer-Verlag
  (2002).

  \bibitem{Br-Bu}
  J. H. Bruinier and M. Bundschuh,
  {\em On Borcherds products associated with
  lattices of prime discriminant}
  to appear in
  Ramanujan J.

  \bibitem{Chen-Yui}
  I. Chen and N. Yui,
  {\em Singular values of Thompson series,}
  in ``Groups, Difference Sets and Monster" (K. T. Arasu {\it et
  al.} Eds.),
  pp. 255-326, de Gruyter, Berlin,
  1995.

  \bibitem{Cohen}
  H. Cohen,
 {\em Sums involving the values at negative integers of $L$-functions
 of quadratic characters},
 Math. Ann. {\bf
 217} (1975),
 271-285.

  \bibitem{C-N}
   J. H. Conway and  S. P. Norton,
   {\em Monstrous Moonshine},
   Bull. London Math. Soc. 11, 308-339,
   1979.

  \bibitem{Cox}
  D. Cox,
  {\em Primes of the form $x^2+ny^2$},
  John Wiley \& Sons,
  1989.

  \bibitem{E-Z}
  M. Eichler and D. Zagier,
  {\em The Theory of Jacobi Forms},
  Progress in Math. {\bf 55},
  Bikh\"auser-Verlag, Boston-Basel-Stuttgart,
  1985.

  \bibitem{G-K-Z}
  B. Gross, W. Kohnen, D. Zagier,
  {\em Heegner points and derivatives of $L$-series II},
  Math. Ann. {\bf 278} (1987),
  497-562.

  \bibitem{Kim}
  C. H. Kim,
  {\em Borcherds products associated with certain Thompson
  series},
  to appear in Compositio Math.,
  e-print arXiv:math.NT/0203076.

  \bibitem{Koike}
  M. Koike,
  {\em On replication formula and Hecke operators},
  Nagoya University
  (preprint).

  \bibitem{Lang}
  S. Lang,
  {\em Algebraic number theory},
  (2nd ed.),
  Springer-Verlag,
  1994.

 \bibitem{Norton}
  S. P. Norton,
  {\em More on moonshine.}
  In: Computational Group Theory, 185-193,
  Academic Press, 1984.
 \bibitem{Zagier94}
 D. Zagier,
 {\em Modular forms and differential operators},
 Proc. Indian Acad. Sci. {\bf 104} (1994),
 57-75.
 \bibitem{Zagier}
 D. Zagier,
 {\em Traces of singular moduli,}
 Max-Planck-Institut f\"ur Mathematik,
 Preprint series {\bf 2000} (8).

\end{thebibliography}
\end{document}